\documentclass[12pt]{amsart}
\usepackage[]{amsmath, amsthm, amsfonts, amssymb, epsfig}

\usepackage{xcolor}
\newcommand\C{{\mathbb C}}

\newcommand\dee{\partial}
\newcommand\Om{\Omega}
\newcommand\Obar{\overline{\Omega}}

\renewcommand\phi{\varphi}

\numberwithin{equation}{section}

\begin{document}

\title[Boundary values]{
A new way to express boundary values in
terms of holomorphic functions on planar Lipschitz domains
}
\author[S.R.~Bell]{Steven R.~Bell}
\author[L.~Lanzani]{Loredana Lanzani}
\author[N.A.~Wagner]{Nathan A. Wagner}

\address{Steven Bell  \hfill\break\indent
Mathematics Department \hfill\break\indent
Purdue University \hfill\break\indent
West Lafayette, IN 47907}
\email{bell@math.purdue.edu}

\address{Loredana Lanzani \hfill\break\indent
Dipartimento di Matematica\hfill\break\indent
Universita' di Bologna
}
\email{loredana.lanzani@gmail.com}

\address{Nathan A. Wagner \hfill\break\indent
 Department of Mathematics \hfill\break\indent
Brown University \hfill\break\indent 151 Thayer St \hfill\break\indent
Providence, RI, 02912 USA}
\email{nathan\_wagner@brown.edu}

\thanks{
L.~ Lanzani is a member of the INdAM group GNAMPA.}
\thanks{N.~A.~Wagner was supported in part by National Science Foundation grant DMS \# 2203272.}


\begin{abstract}
We decompose $p$ - integrable functions on
the boundary of a simply connected Lipschitz domain $\Omega \subset \mathbb C$
into the sum of the boundary values of two, uniquely determined holomorphic functions, where
one is holomorphic in $\Omega$  while the other is holomorphic in $\mathbb C\setminus \overline{\Omega}$ and
vanishes at infinity. This decomposition has been
described previously for smooth functions on the boundary of
a smooth domain \cite{Be}. Uniqueness of the decomposition is elementary in the smooth
case, but extending it to the $L^p$
setting relies upon a classical albeit little-known regularity theorem for the holomorphic Hardy space $h^p(b\Om)$ of planar domains for which we provide a new proof that is valid also in higher dimensions. An immediate consequence of our result will be a new
characterization of the kernel of the Cauchy transform acting on $L^p(b\Omega)$. These results give a new perspective on the classical
Dirichlet problem for harmonic functions and the Poisson formula even in the case of the
disc. Further applications are presented along with 
directions for future work.
\end{abstract}

\maketitle

\theoremstyle{plain}

\newtheorem {thm}{Theorem}[section]
\newtheorem {lem}[thm]{Lemma}
\newtheorem{cor}[thm]{Corollary}

\hyphenation{bi-hol-o-mor-phic}
\hyphenation{hol-o-mor-phic}

\section{Introduction}
\label{sec1}
If $\Om$ is a simply connected domain in the plane bounded
by a $C^\infty$ smooth Jordan curve, then there are various
ways to express a $C^\infty$ smooth function $u$ given on
the boundary in terms of two holomorphic functions on 
$\Om$ that extend $C^\infty$ smoothly to the boundary.
The most natural way is to require that
$$u=f+\overline{F}$$
on the boundary. By choosing a point $a$ in $\Om$ and
further requiring that $F(a)=0$, the decomposition is
unique and yields the solution to the classical Dirichlet
problem for harmonic functions with boundary data $u$, \cite{Be}.

If the boundary
curve of $\Om$ is parametrized in the counterclockwise sense
by $z(t)$, $a\le t\le b$, the complex unit tangent function
$T(z)$ on the boundary is defined via $T(z(t))=z'(t)/|z'(t)|$.
For each point $z$ in the boundary, $T(z)$ is a complex number
of unit modulus pointing in the direction of the tangent
vector at $z$ in the direction of the standard orientation.
Note that $i\,T(z)$ represents the direction of the inward
pointing normal vector to the boundary at $z$. Here and throughout
we will work with arc-length measure $\sigma$ without explicitly denoting it in
the function spaces; hence $L^p(b\Om)$ stands for $L^p(b\Om, d\sigma)$ etc.

A second
way to express the function $u$ on the boundary in terms
of holomorphic functions on $\Om$ is
$$u=g+\overline{G\,T}.$$
This gives an orthogonal decomposition of $u$ in terms
of a holomorphic function $g$ and another function
$\overline{G\,T}$ that is orthogonal to the Hardy space
$h^2(b\Omega)$ in the $L^2(b\Om)$ inner product. The function $g$ is equal
to the Szeg\H o projection of $u$ and
$\overline{G\,T}$ is the projection of $u$ onto the
orthogonal complement of $h^2(b\Om)$ in $L^2(b\Om)$. The functions $g$ and $G$ are uniquely
determined.

The relationships between the functions in the decompositions
described thus far are fascinating and have been studied
extensively. In this paper, we study a third way to express the
boundary function $u$, namely
 $$
 u= h+H$$ where $h$ is holomorphic in $\Om$, while
$H$ is holomorphic on $\C\setminus\Obar$, the ``outside of $\Om$,''
and tends to zero at infinity. We call this the $h+H$ decomposition from now on.
This decomposition can be viewed as a direct consequence of the ``jump
theorem''  that states that the jump of the Cauchy transform of $u$ across a boundary point
 $z$ is the value of $u$ at $z$ for almost every $z\in b\Om$. This classical fact is a direct consequence of the Plemelj
formula. To better understand it, we will look at it from several
points of view, including a way that yields a higher order Plemelj
formula in the smooth setting.
We take special care in extending the $h+H$ decomposition to the setting of Lipschitz domains: this is because the first
 term, $h$, turns out to be the Cauchy transform of the data $u$, which is known
 to be bounded on $L^p(b\Om)$ for a large class of non-smooth domains that include the Lipschitz domains, see \cite{D}; and the second term, $H$ is
given by a so-called ``solid Cauchy transform''
formula in \cite[p.~87]{Be}. 
 Historically, the Lipschitz category has brought forward a paradigm shift in the analysis of the Cauchy transform with the seminal works \cite{Ca} and \cite{CMM}, which in turn have paved the way to a fully developed arsenal of tools that allow for optimal results, see e.g., \cite[Theorem 2.1]{LS}. Furthermore, in applied mathematics 
 non-smooth domains often provide more desirable settings than the $C^\infty$ category, see \cite{HLLS1} and \cite{HLLS2}.

The main result in this paper is the following decomposition theorem, where we use the notation $\vartheta(\Om)$ for the set of functions holomorphic on a given domain $\Om$.

\begin{thm}\label{T:1}
Suppose $\Om$ is a bounded Jordan domain with Lipschitz boundary.
\begin{itemize}
\item[{\tt i.}] For $1< p < \infty$, let $h^p(b\Om)$ and $h^p_-(b\Om)$ denote the boundary values of the Hardy spaces for $\Om$ and for $\mathbb C\setminus \Obar$, respectively. 
If $u\in L^p(b\Om)$ then $u$ can be expressed uniquely as $h+H$ a.e. $b\Om$ where
$h\in h^p(b\Om)$ and $H\in h^p_-(b\Om)$.
\vskip0.1in

\item[{\tt ii.}] If, moreover,  $u \in C^{0,\alpha}(b \Om)$ for some $0<\alpha < 1$, then $u$ can be expressed uniquely as $h+H$ where $h\in \vartheta(\Om)\cap\, C^{0,\alpha}(\overline{\Om})$, while $H$ is in $\vartheta(\mathbb C\setminus\overline{\Om})\cap\, C^{0,\alpha}(\mathbb C\setminus{\Om})$ and vanishes at infinity.
\vskip0.1in

\item[{\tt iii.}] Suppose additionally that $\Om$ is of class $C^{k,\alpha}$ for some $k \geq 1$ and $0<\alpha < 1.$ If $u \in C^{k,\alpha}(b \Om)$, then $u$ can be expressed uniquely as $h+H$ where
$h\in \vartheta(\Om)\cap\, C^{k,\alpha}(\overline{\Om})$, while
$H$ is in $\vartheta(\mathbb C\setminus\overline{\Om})\cap\, C^{k,\alpha}(\mathbb C\setminus{\Om})$ and vanishes at infinity.
\vskip0.1in

\item[{\tt iv.}] Suppose additionally that $\Om$ is of class $C^\infty$. If $u$ is in $C^\infty (b\Om)$,
then $u$ can be expressed uniquely as $h+H$ where
$h\in \vartheta(\Om)\cap\, C^{\infty}(\overline{\Om})$, while
$H$ is in $\vartheta(\mathbb C\setminus\overline{\Om})\cap\, C^{\infty}(\mathbb C\setminus{\Om})$ and vanishes at infinity.

\end{itemize}
In any of the cases {\tt i.} to {\tt iv.}, we have
\begin{equation}\label{E:h}
h(z) = \mathcal C u(z)
\end{equation}
where $\mathcal C$ is the Cauchy transform associated to $\Om$. 
\end{thm}
It follows that
\begin{equation}\label{E:KerCauchy}
\mathcal C u =0 \quad\iff\quad u = H
\end{equation}
where $H$ has the above-stated properties.

Theorem \ref{T:1} admits the following corollary, which gives a Banach direct sum decomposition of $L^p(b\Om)$ in terms of internal and external Hardy spaces. 

\begin{cor}
Suppose $\Om$ is a bounded Jordan domain with Lipschitz boundary and $1<p<\infty$. Then the following Banach space isomorphism holds:

\begin{equation} L^p(b \Om) \simeq h^p(b\Om) \oplus h^p_{-}(b \Om) \end{equation}

\begin{proof}
Apply Theorem \ref{T:1}, part {\tt i.}, together with the $L^p$ estimate for the Cauchy transform $\mathcal{C}$ obtained in \cite{CMM}. 
\end{proof}
\end{cor}

Conclusion
{\tt iv.} and the existence part of conclusion {\tt ii.} were studied
in \cite[p.~87]{Be} in the case when $\Om$ is $C^\infty$-smooth: we briefly recall the argument given there. Given $u\in C^\infty(b\Om)$ the decomposition: $u=h+H$ is based on the Cauchy-Pompeiu integral formula
$$\Phi(z)=\frac{1}{2\pi i}\int_{b\Om}\frac{\Phi(w)}{w-z}\ dw+
\frac{1}{2\pi i}
\iint_{\Om}\frac{\frac{\dee \Phi}{\dee\bar w}(w)}{w-z}\ dw\wedge d\bar w,\quad z\in \Om$$
which holds for a function $\Phi$ in $C^\infty(\Obar)$, and
the observation that there exist functions $\phi_n$ in $C^\infty(\Obar)$
that vanish on the boundary such that $\dee\phi_n/\dee\bar w$
and $\dee \Phi/\dee\bar w$ agree to (arbitrarily high) order $n$
on the boundary (see Lemma~2.1 in \cite[p.~6]{Be} for a proof
and more details about how we are about to use this lemma).

Thus, if we extend
our given $u$ in $C^\infty(b\Om)$ to (any) $U$ in $C^\infty(\Obar)$
and then pick a function $\phi_n$ in $C^\infty(\Obar)$ that vanishes
on the boundary such that
$\Psi_n:=\dee U/\dee\bar w - \dee\phi_n/\dee\bar w$ vanishes to
order $n$ on the boundary, we can apply the Cauchy-Pompeiu
formula to $\Phi := U-\phi_n$ to see that
\begin{equation}\label{E:1}
U(z)-\phi_n(z)=\frac{1}{2\pi i}\int_{b\Om}\frac{u(w)}{w-z}\ dw+
\frac{1}{2\pi i}
\iint_{\Om}\frac{\Psi_n(w)}{w-z}\ dw\wedge d\bar w,\quad z\in \Om
\end{equation}
and we deduce that the function
defined by the second integral in \eqref{E:1} defines a
holomorphic function $H$ on the outside of $\Om$ that vanishes at
infinity and that extends $C^n$ smoothly to the boundary
of $\Om$ from the outside. Indeed, we may extend $\Psi_n$ by
zero outside $\Om$ to obtain a $C^n$ function on $\C$ and
then change variables in the integral on $\Om$ to write it as
$$\frac{1}{2\pi i}
\iint_{\C}\frac{\Psi_n(z-w)}{w}\ dw\wedge d\bar w$$
to see that we can differentiate in $z$ under the integral to
conclude that the above defines a $C^n$ function on $\C$.
Thus, we can read off that the boundary integral in \eqref{E:1} also extends $C^n$ smoothly to
the boundary from the inside, and we call its boundary value 
$\mathcal C u$, the Cauchy transform of $u$. That is, the formula
\eqref{E:1} extends to $\Obar$ and on the boundary $b\Om$ it takes the form

\begin{equation}
u(z)=\mathcal Cu(z)+
\frac{1}{2\pi i}
\iint_{\Om}\frac{\Psi_n(w)}{w-z}\ dw\wedge d\bar w,\quad z\in b\Om.
\end{equation}

Thus $u=h+H$ where 

$$
h(z)  :=\mathcal C u (z)\quad \text{and}\quad 
H(z) := \frac{1}{2\pi i}
\iint_{\Om}\frac{\Psi_n(w)}{w-z}\ dw\wedge d\bar w.
$$
 
Both functions are holomorphic on their respective sides of the
boundary and extend to be $C^n$ up to the boundary from their
respective sides.
Note that integration by parts reveals that
$$\iint_{\Om}\frac{\frac{\dee\phi_n}{\dee\bar w}(w)}{w-z}\ dw\wedge d\bar w\ =\ 0,\quad z\in \mathbb C\setminus \Om\, .$$

 It follows that
$$H(z)=
\frac{1}{2\pi i}
\iint_{\Om}\frac{\frac{\dee U}{\dee\bar w}(w)}{w-z}\ dw\wedge d\bar w,\quad z\in \mathbb C\setminus \Obar .$$
That $H$ is independent of the choice
of the extension of $u$ to $\Obar$ is clear because, if $U_1$
and $U_2$ were two such extensions, then $U_1-U_2$ would vanish on
the boundary and the integration by parts argument applied to
$\dee U_1/\bar w- \dee U_2/\bar w$ would yield the independence
of the definition of $H$ on the choice of $U$. Since $n$ was
arbitrary, it follows that $\mathcal C u$ and $H$ are $C^\infty$
smooth up to the boundary from their respective sides. The last
expression for $H$ can be further reduced by noting that
$\frac{\dee U}{\dee\bar w}(w)/(w-z)=\frac{\dee}{\dee\bar w}
\left[U(w)/(z-w)\right]$ and applying integration by parts to obtain
$$H(z)=
-\frac{1}{2\pi i}
\int_{b\Om}\frac{u(w)}{w-z}\ dw,$$
and thus $H$ can be viewed as ``minus the Cauchy transform of $u$
on the other side of the boundary'' and we see that this approach
to the subject yields a way to understand the Plemelj formula
from a $C^\infty$ perspective.

As observed in \cite[p.~87]{Be}, for smooth $\Om$ the decomposition
$u=h+H$ in item {\tt i.} is unique because,
if we had two such expressions $u=h_1+H_1=h_2+H_2$, then
$h_1-h_2=H_2-H_1$ on the boundary, and so an extension of
$h:=h_1-h_2$ to the complex plane as an entire function is obtained by
defining the extension to be $H:=H_2-H_1$ on the outside and
noting that the functions glue together along the boundary
and so are holomorphic across the boundary via a standard
Morera's theorem argument. The vanishing of the holomorphic functions
on the exterior of $\Om$ at infinity shows that the extension is a
bounded entire function, and so Liouville's theorem implies that
$H$ must be the zero function. This concludes the proof of {\tt i.} for smooth 
$\Om$ as given in \cite{Be}. 

Concerning the proof of the existence part of conclusion {\tt ii.} for smooth $\Om$: as noted in \cite{Be}, the boundedness of the Cauchy
transform as an operator from $L^p(b\Om)$ to $h^p(b\Om)$ for $1<p<\infty$, and the density in $L^p(b\Om)$ of
$C^\infty(\Obar)$, reveals that, if $u$ is merely in $L^p(b\Om)$,
then $h=\mathcal C u$ is in $h^p(b\Om)$. The techniques used
in \cite{Be} also yield that $H$ on the boundary is an
$L^p$ limit of the boundary values of holomorphic functions
on the outside of $\Om$ that are $C^\infty$ smooth up to
$b\Om$ and that vanish at infinity. It can be seen that
the operator taking $u$ to $H$ is also bounded in $L^p(b\Om)$
and that $H$ is in $h^p_-(b\Om)$ that is, $H$ is given as the $L^p$ boundary values of a
function that is holomorphic on the outside of $\Om$ and
vanishes at infinity. This version of Hardy space on the
outside of $\Om$ can be understood by mapping the outside
of $\Om$ to a bounded domain via inversion and
then using standard results about convergence in the Hardy
space of a bounded simply connected domain with smooth
boundary on the subspace of functions that vanish at the
point that is the image of the point at infinity under the
inversion.

The uniqueness of the $h+H$ decomposition for $u$ merely in $L^p(b\Om)$ was well-known for $p=2$ in the special case when $\Om = \mathbb D$ (the unit disc): here uniqueness of the $h+H$ decomposition for $u$ in $L^2$ can be derived from a classical result of Hardy and
Littlewood that states that a complex antiderivative of a
function in $h^2(b\mathbb D, d\theta)$ is H\"older-$1/2$ continuous
on $\mathbb D$ and extends to be H\"older-$1/2$ continuous on $\overline{\mathbb D}$. (Another way to deduce that the extension is entire in the
uniqueness argument follows from results of Koosis in
\cite[pp.~87-88]{K} about Schwarz reflection in the $L^2$
setting on the disc.)
\medskip

Theorem \ref{T:1} {\tt iv.} has the following corollary.
\begin{thm}\label{T:2}
Suppose $\Om$ is a smooth Jordan domain. 
\begin{itemize}
\item[{\tt i.}] Let $U$ be in $C^\infty(\overline{\Om})$. Then there is $\varphi$ in $C^\infty(\overline{\Om})$ that vanishes on $b\Om$ and such that $\displaystyle{\frac{\dee\varphi}{\dee\bar z}}$ agrees with $U$ on $b\Om$ up to infinite order.
\medskip

 \item[{\tt ii.}] Let $V$ be in $C^\infty(\mathbb C)$. Then there is $\Phi$ in $C^\infty(\mathbb C)$ that vanishes on $b\Om$  and such that $\displaystyle{\frac{\dee\Phi}{\dee\bar z}}$ agrees with $V$ on $b\Om$ up to infinite order.
 \end{itemize}
 \end{thm}
 Indeed, by Theorem \ref{T:1} we have expressed
the boundary values $u$ of $U$ as $h+H$.
Let $\widetilde H$ denote a $C^\infty$ extension of $H$ to
$\C$ that agrees with $H$ outside $\Om$ and let $\tilde h$ denote
a $C^\infty$ extension of $h$ to $\C$ that agrees with $h$ inside
$\Om$. Note that
$\dee \tilde h/\dee\bar w$ and
$\dee \widetilde H/\dee\bar w$
vanish to infinite order on the boundary since they are
identically zero on one side of the boundary.
Now $\varphi :=U-\tilde h - \widetilde H$ restricted to $\Obar$
satisfies our requirements. Such a $\varphi$ could also be
constructed via PDE techniques by using the construction
of functions satisfying arbitrarily high degree of vanishing
conditions and standard arguments based on Borel's lemma
(see \cite[p.~98]{GG}). Another approach is via a Mittag-Leffler
argument like the ones described in \cite[p.~201]{Be1}
It is also clear that local versions of the theorem for
$C^\infty$ functions $v$ on $C^\infty$ curves $\gamma$ also
follow.

If we started with a function $V$
that is $C^\infty$ on $\C$, we could modify
$V$ by multiplying it by a $C^\infty$ cut-off function that
is compactly supported and equal to one on a big disc containing
$\Obar$, then solve a $\dee$-bar problem
$\dee U/\dee\bar w = v$ via
$$U(z)=
\frac{1}{2\pi i}
\iint_{\C}\frac{v(w)}{w-z}\ dw\wedge d\bar w$$
and then construct $\Phi$ in   
$C^\infty(\C)$ using such $U$.
\bigskip

It is worthwhile pointing out that the uniqueness of the $h+H$ decomposition in $L^p(b\Om)$ relies on 
the fact that for any $p
\geq 1$, all members of $h^p(b\Om)$ turn out to be complex derivatives of functions continuous up to the boundary of $\Om$, and on an analogous statement for $h^p_-(b\Om)$. Namely

\begin{thm}\label{T:3}
Suppose $\Om$ is a bounded Jordan domain with Lipschitz boundary

and let $1 \leq p<\infty$.

\begin{itemize}
\item[{\tt i.}] For any $h$ in $h^p(b\Om)$ there is $H \in \vartheta(\Om) \cap C(\Obar)$ such that
$$
H'(z) = h (z) \quad \text{for any}\quad z\in \Om\, ,
$$
and $H \lvert_{b \Om}$ is absolutely continuous with respect to arc length on $b \Om$.
\item[{\tt ii.}] For any $h$ in $h^p_-(b\Om)$ there is $H \in \vartheta(\mathbb C\setminus \Obar)\cap C(\mathbb C\setminus \Om)$ that vanishes at infinity such that
$$
H'(z) = h (z) \quad \text{for any}\quad z\in \mathbb C\setminus\Obar\, ,
$$
and $H \lvert_{b \Om}$ is absolutely continuous with respect to arc length on $b \Om.$
\end{itemize}
\end{thm}
This is a partial restatement of a classical theorem of Privalov \cite[Sect. 7.3]{P}. Here we present a new proof which, unlike the original argument in \cite{P}, does not rely on conformal mapping and is thus amenable to be extended to higher dimensions. We remark that this result holds for the broadest (Banach space) class $h^1(b\Om)$ (that is for $p=1$), unlike the decomposition result Theorem \ref{T:1}. This is because the decomposition result naturally relies on the boundedness of the Cauchy transform for $p>1$, while the endpoint estimate at $L^1$ is well-known to fail.  

\medskip

A final observation is in order. Since the $h+H$ decomposition for $u\in L^2(b\mathbb D)$ that is, for the case when $\Om$ is the disc, is a trivial consequence of the Fourier series decomposition of $u$, it is legitimate to ask whether a proof for general $\Om$ can be obtained via a conformal mapping argument. It turns out that a conformal map-based scheme would employ two such maps: a conformal map $f: \mathbb D\to \Om$ along with a conformal map $f^*: \mathbb C\setminus \overline{\mathbb D}\to\mathbb C\setminus \overline{\Om}$, and one would need to have
$$
f(\zeta) = f^*(\zeta)\qquad \text{for}\quad |\zeta|=1.
$$
But if this were indeed the case, a theorem of Besicovich \cite[Theorem 3.4]{Younsi} would force $f$ and $f^*$ to be (the same) constant, giving a contradiction.
\bigskip

This paper is organized as follows. In Section \ref{sec2} we prove Theorem \ref{T:3} and the uniqueness of the $h+H$ decomposition in $L^p(b\Om)$  under the additional assumption that $\Om$ is smooth (of class $C^\infty$); we present the treatment of smooth domains in a separate section because the proofs are of independent interest and because Lipschitz domains can be exhausted by families of smooth sub- and super-domains, a well-known fact which is relevant here and is recalled in Section \ref{sec3}.

 Theorem \ref{T:3} and our main result, Theorem \ref{T:1}, are proved in Section \ref{sec3}, where we also review the main features of the holomorphic Hardy spaces for Lipschitz $\Om$ and for (the interior of) its complement. Finally, in Section
\ref{sec4} we collect a few applications of the $h+H$ decomposition yielding a new representation for the solution of the classical Dirichlet problem for harmonic functions on the disc (Theorem \ref{T:4} and Remark 2), and a new pseudo-local property of the Szeg\H o projection for $C^\infty$-smooth $\Om$ (Theorem \ref{T:5}).
\medskip

\noindent{\bf Acknowledgement.} We are deeply grateful to Dmitry Khavinson for bringing to our attention a theorem of Privalov \cite{P} which was unknown to us.

\bigskip

\section
{Smooth domains: antiderivatives
 of  Hardy space functions and uniqueness of the $h+H$ decomposition in $L^p(b\Om)$}
\label{sec2}
\subsection{Proof of Theorem \ref{T:3} for smooth $\Om$}

\label{Ss:2.1}

We let $\gamma$ denote the boundary of $\Om$
parametrized in the counterclockwise sense via $z(t)$, 
$a\le t\le b$.
 
We may
assume that $z_0=z(a)$ is a point in the boundary of $\Om$
where $h$ has a non-tangential boundary limit. For small
$\epsilon>0$, let $\gamma_\epsilon$ denote the ``parallel''
curve to $\gamma$ parametrized by
$z_\epsilon(t)=z(t)+
i\epsilon T(z(t))$
gotten by moving a distance $\epsilon$
along the inward pointing normal to the boundary curve of
$\Om$ from $z(t)$ as $z(t)$ traces out the boundary. It is
well know that, if $\epsilon>0$ is sufficiently small, then
$\epsilon$ is equal to the distance from $z_\epsilon(t)$ to
the boundary, and $z(t)$ is the closest point in the boundary
to $z_\epsilon(t)$. We will assume henceforth that our future
epsilons are in this range.

Before we proceed further, we define a closed contour in
$\Om$ that will be useful. If $0<\hat\epsilon<\epsilon$, we
let
$\mathcal R(\hat\epsilon,\epsilon,\tau)$ denote the
curvelinear rectangle that
starts at $z_0+i\hat\epsilon T(z_0)$ and
follows $\gamma_{\hat\epsilon}$ as the parameter $t$ ranges from
$a$ to $\tau$. When the inward normal to $z(\tau)$ is reached,
the curve turns left and follows the normal to
$z(\tau)+i\epsilon T(z(\tau))$. It then turns left again and
follows $\gamma_{\epsilon}$ backwards to
$z_0+i\epsilon T(z_0)$, and then finally follows the inward
pointing normal at $z_0$ back to the starting point at
$z_0+i\hat\epsilon T(z_0)$. The Cauchy theorem implies that
the integral of $h(z)$ around 
$\mathcal R(\hat\epsilon,\epsilon,\tau)$ is zero.

First prove the result for $p=2$. The $L^2$ convergence of $h$ along the curves
$\gamma_\epsilon$ to the boundary values of $h$ along
$\gamma$ implies that the integrals of $h$ along the curved
parts of $\mathcal R(\hat\epsilon,\epsilon,\tau)$ roughly
cancel, and their sum tends to zero as $\epsilon$ tends
to zero. We note further that, since $h$ has a limit along the inward
pointing normal at $z_0$, that the integral of $h$ along a
segment of the inward pointing normal of length
$\epsilon-\hat\epsilon$ goes to zero as $\epsilon$ tends to
zero. Hence, because the integral of $h(z)$
around $\mathcal R(\hat\epsilon,\epsilon,\tau)$ is zero, and
because the integrals along the curved parts of the contour
roughly cancel, integrals of $h$ along a segment of the inward
pointing normal to $z(\tau)$ of length $\epsilon-\hat\epsilon$
tend uniformly in $\tau$ to zero as $\epsilon$
tends to zero.

We can now conclude that that the integral of $h(z)$ around
$\mathcal R(0,\epsilon,\tau)$ is well defined if we express
the part of the integral along the inward pointing normal
at $z(\tau)$ as a limit as $\hat\epsilon$ tends to zero.
Now, the vanishing of the integral of $h(z)$ around
$\mathcal R(\hat\epsilon,\epsilon,\tau)$
implies that the integral of $h(z)$ around
$\mathcal R(0,\epsilon,\tau)$ is also equal to zero.

We now define a holomorphic function $H(z)$ on $\Om$, complete
with continuous boundary values, whose complex derivative
is equal to $h(z)$ in $\Om$. We start by defining $H(z)$ on
the boundary via
$$H(z(\tau))=\int_a^\tau h(z(t))z'(t)\ dt$$
Because the boundary values of $h$ are in $L^2(b\Om)$, the
Cauchy-Schwarz inequality can be used to easily show that $H$ is
continuous on the boundary. Indeed, applying Cauchy-Schwarz to
the integral of $h$ along a curve segment $\sigma$ in the boundary
and parameterizing the segment in terms of arc-length yields that
$$\left|\int_\sigma h\ dz\right|\le \|h\|\ell(\sigma)^{1/2}$$
where $\|h\|$ is the $L^2$ Hardy space norm of $h$ on the whole
boundary and $\ell(\sigma)$ is the length of $\sigma$.
Note that $H(z_0)=0$. The fact that $\int_\gamma h\, dz=0$
(by the Cauchy theorem) guarantees that the value of $H(z(t))$ as
$t$ tends to the endpoint $t=b$ is also zero. Moreover, $H$ is in fact absolutely continuous on the boundary by virtue of its definition as the integral of an $L^1$ function.

Next, we define $H(z)$ at a point inside $\Om$ as follows.
Let $\eta_z$ denote a curve that starts at $z_0$ and follows
the inward pointing normal to the boundary at $z_0$ the short
distance $\epsilon$, and then follows any curve in
$\Om$ from the terminus of the short interval to $z$. It
is easy to see that
$$H(z)=\int_{\eta_z}h(w)\ dw$$
defines an analytic antiderivative of $h$ on $\Om$ and
that the limit of $H(z)$ along the inward pointing normal
to the boundary at $z_0$ is equal to zero. Note also that
the definition of $H$ is independent of the choice of $\epsilon$.

We now prove that $H_\epsilon(t):=H(z_\epsilon(t))$ tends uniformly to
$H_0(t):=H(z(t))$ as $\epsilon$ tends to zero through positive
values. We have shown that the integral of $h$ around
$\mathcal R(0,\epsilon,\tau)$ is zero and that the integral of $h$
on the line
$L(\epsilon,\tau)$ from $z_\epsilon(\tau)$ to $z(\tau)$ is well
defined and tends to
zero as $\epsilon$ tends to zero and this convergence is uniform
in $\tau$ because of the $L^2$ convergence of $h(z_\epsilon(t))$ to
$h(z(t))$ on the whole boundary, and the fact that the integral
along the line segment along the inward pointing normal to $z_0$
tends to zero. We may now conclude that
$$ H(z(\tau))-H(z_\epsilon(\tau))
=\int_{\mathcal R(0,\epsilon,\tau)}h\ dz+
\int_{L(\epsilon,\tau)}h\ dz=
\int_{L(\epsilon,\tau)}h\ dz,$$
and the uniform convergence follows. Since the convergence is
uniform, and since the boundary values are continuous, we conclude
that $H$, as we have defined it on $\Obar$, is continuous there.
\medskip

The proof for $p>1$ is similar, 
because functions in 
$h^p(b\Om)$ have $L^p$ boundary values and H\"older's inequality
can be used to replace Cauchy-Schwarz inequality in the calculation
above to see, in case $p>1$, that the boundary values of the
antiderivative are continuous, since
$$\left|\int_\sigma h\ dz\right|\le \|h\|\ell(\sigma)^{1/q}$$
where $\|h\|$ is the $L^p$ Hardy space norm of $h$.
If $p=1$, this integral goes to zero as the length $\ell(\sigma)$
goes to zero by measure theory principles.

This method of proof  carries over to $h^p_-(b\Om)$ (the Hardy space for the complement
of $\Obar$).
\subsection{Uniqueness of the $h + H$ decomposition in $L^p(b\Om)$ for smooth $\Om$}\label{Ss:2.2}

Suppose $u=h_1+H_1=h_2+H_2$. Then the boundary values of $h=h_1-h_2$
and $H=H_2-H_1$ agree almost everywhere. A point
$z_0$ in the boundary can be chosen so that both $h$ and $H$ have
non-tangential boundary limits at $z_0$ from their respective sides.
Antiderivatives of the respective functions on both sides can be defined
via path integrals in each side that start at the point $z_0$
and follow the inward pointing normals a short distance before
heading to a point in the domain. When the construction of the
antiderivatives of $h$ and $H$ are carried out, and the
continuity up to the boundary is noted, the resulting
functions are seen to agree on the boundary. Consequently, an
entire function is obtained. Its derivative reveals that $h$
extends via $H$ to be entire and tends to zero at infinity, and so it is
identically zero. The uniqueness of the decomposition follows.
We expect some of these details of this version of a proof to be
useful in our future work as we generalize the decomposition to
more general functions on more general domains using tools from
\cite{LS}.

{\bf Remark 1.} The results of this section can
be localized to yield a rather simple proof of a removable
singularity theorem that states that if a function is holomorphic
on both sides of a smooth curve and has $L^2$ boundary limits from
both sides that agree along the curve, then the function extends
holomorphically across the curve. Koosis proved this result
in the disc using the Schwarz function $S(z)=1/z$ of the disc,
but we suspect it must be a known result for any smoothly bounded $\Om$.

\section{Proof of Theorem \ref{T:1}}
\label{sec3}
 
\subsection{Background on Lipschitz Domains}\label{ssec3.1}
 We say a bounded domain $\Omega$ is Lipschitz if there exists a finite collection of coordinate rectangles $\{R_j\}_{j=1}^{N}$ so that $R_j \cap b\Omega$ is the graph of a real-valued Lipschitz function for each $j.$ Precisely, this means that for each $j$, there exists a Lipschitz function $\varphi_j:\mathbb{R} \rightarrow \mathbb{R}$ and angle $\theta_j \in [0,2 \pi)$ so that $$R_j \cap \Omega= R_j \cap \{e^{i \theta_j}(x+ iy): y> \varphi_j(x)\}.$$

For such domains there exists a Lipschitz diffeomorphism from $b\Omega$ to $b \widetilde{\Omega}$, where $\widetilde{\Omega}$ is a $C^\infty$ smooth domain. Let us additionally suppose that $\Omega$ is simply connected and bounded by a Jordan curve. By the above fact, we see that there exists a bi-Lipschitz parametrization $z(t):[a,b] \rightarrow b\Omega$ of said boundary curve. The derivative $z'(t)$ exists for almost every $t \in [a,b]$ and as usual unit tangent vector $T(z(t))$ is given by $\frac{z'(t)}{|z'(t)|}$. We will again assume that the parametrization is taken in the counterclockwise direction. 

In particular, there exists a family of $C^\infty$ smooth domains $\{\Omega_\epsilon\}_{\epsilon>0}$ that approximate $\Omega$ non-tangentially from the inside (that is, $\{\Om_\epsilon\}_\epsilon$ is an exhaustion of $\Om$) in a quantitatively precise sense, and likewise from the exterior, see \cite{VerchotaThesis} for the full details of the construction (originally due to Necas). Here, we point out that the family $\{\Omega_\epsilon\}_\epsilon$ depends continuously on the parameter $\epsilon$ (and so we henceforth refer to it as ``a continuous family'') in the sense that any parametrization $t\mapsto \gamma_\epsilon(t) $ for the boundary $b\Om_\epsilon$ is a continuous function of $(t, \epsilon)$. Note that such approximation is stated in \cite{VerchotaThesis,Verchota, L2000} only in terms of a sequence $\{\Omega_j\}$, but the proof in fact gives a continuous family from which a discrete one is extracted. The boundaries of the approximating domains are given by the level sets of a family of functions $\chi_\epsilon$ which are in turn obtained by convolving the characteristic function of $\Omega$ with a family $\{\psi_\epsilon\}$ of smooth approximations to the identity (which again depends continuously on $\epsilon$).

The non-tangential approach regions for $\Omega$ can in fact be taken to be truncated cones. Given any point $p \in b \Omega,$ we use $\Gamma(p)$ to denote a truncated double cone (with two non-empty components) with vertex at $p$. Let $\Gamma_i(p)$ denote the component of the cone $\Gamma(p)$ that intersects $D$, and $\Gamma_e(p)$ the component that intersects $\mathbb{C} \setminus \overline{D}$. The apertures of these cones can be bounded above and below for a Lipschitz domain to guarantee each component lies totally inside or totally outside $D$. In other words, there exists a family of such double cones $\{\Gamma(p)\}_{p \in b \Omega}$ satisfying the following properties: There exist truncated double cones $\alpha, \beta$, centered at the origin and whose axis is the $y$-axis, so that if $\alpha_i, \beta_i$ denote the components of $\alpha, \beta$ above the $x$-axis and $\alpha_e, \beta_e$ the components below the $x$-axis, then we have, for each $1 \leq j \leq N$ and $p \in R_j \cap \partial \Omega$:

$$ e^{i \theta_j} \alpha_i+p \subset \Gamma_i(p) \subset \overline{\Gamma_i(p)} \setminus \{p\} \subset e^{i \theta_j} \beta_i+p \subset D \cap R_j;$$ and 

$$ e^{i \theta_j} \alpha_e+p \subset \Gamma_e(p) \subset \overline{\Gamma_e(p)} \setminus \{p\} \subset e^{i \theta_j} \beta_e+p \subset (\mathbb{C} \setminus \overline{D}) \cap R_j.$$

Therefore, we have $\Gamma_i(p)= \Gamma(p) \cap D$ and $\Gamma_e(p)= \Gamma(p) \cap \mathbb{C} \setminus \overline{D}$, so importantly we see that these interior and exterior non-tangential approach regions are convex (they are components of a truncated cone). 

For future reference we record below the details concerning the interior and exterior approximating domains:

There exist two families of $C^\infty$-smooth domains $\Omega_\epsilon \subset \Omega$ and $\Omega^\epsilon \supset \Omega$ approximating $\Omega$ in the following sense (for sufficiently small $\epsilon>0$):
\begin{enumerate}
\item There exists a family of of Lipschitz diffeomorphisms
$$\Lambda_\epsilon: b\Omega \rightarrow b \Omega_\epsilon\qquad \text{and}\qquad \Lambda^\epsilon: b \Omega \rightarrow b \Omega^\epsilon$$ such that
$$ \sup \{|\Lambda^\epsilon(p)-p|+|\Lambda_\epsilon(p)-p|: p \in b \Omega \} \rightarrow 0,\quad  \epsilon \rightarrow 0.$$  
\vskip0.1in

Moreover, one has $\Lambda_\epsilon(p) \in \Gamma_i(p)$, $\Lambda^\epsilon(p) \in \Gamma_e(p)$ for each $\epsilon.$ 
\vskip0.1in
\item We may assume the original family $\{R_j\}$ also forms a  family of coordinate rectangles for $b\Omega_\epsilon, b \Omega^\epsilon$ for each $\epsilon$. Moreover, for every such rectangle $R_j$, if $\varphi, \varphi_\epsilon, \varphi^\epsilon$ denote the Lipschitz functions whose
graphs describe the boundaries of $ \Omega$, $ \Omega_\epsilon$, and $\Omega^\epsilon$ respectively, in $R_j$, then $\|\nabla \varphi_\epsilon\|_{L^\infty}, \|\nabla \varphi^\epsilon\|_{L^\infty}  \leq  \|\nabla \varphi\|_{L^\infty}$ and $\nabla \varphi_\epsilon, \nabla \varphi^\epsilon \rightarrow \nabla \varphi$ a.e. and in every $L^q$, $1<q<\infty.$
\vskip0.1in
\item There exist positive functions $w_\epsilon, w^\epsilon: b \Omega \rightarrow \mathbb{R}_{+}$, bounded away from zero and infinity uniformly in $\epsilon$, such that for
any measurable set $F \subset b\Omega$, $\int\limits_{F} w_\epsilon \, d \sigma= \int\limits_{\Lambda_\epsilon(F)} d \sigma_\epsilon$ and $\int\limits_{F} w^\epsilon \, d \sigma= \int\limits_{\Lambda^\epsilon(F)} d \sigma^\epsilon$, where $d \sigma_\epsilon$, $d \sigma^\epsilon$ denote arc-length measure on $b \Omega_\epsilon$, $b \Omega^\epsilon$, respectively. The following
change of variables formulas hold:
\begin{equation}
\int\limits_{b \Omega} f_\epsilon(\Lambda_\epsilon(p)) w_\epsilon(p) \, d\sigma(p) = \int\limits_{b \Omega_\epsilon} f_\epsilon(p_\epsilon) \, d \sigma_\epsilon(p_\epsilon) \label{COV};
\end{equation}

\begin{equation}
\int\limits_{b \Omega} f^\epsilon(\Lambda^\epsilon(p)) w^\epsilon(p) \, d\sigma(p) = \int\limits_{b \Omega^\epsilon} f^\epsilon(p^\epsilon) \, d \sigma^\epsilon(p^\epsilon).
\end{equation}
where $f_\epsilon, f^\epsilon$ is any measurable
function on $b \Omega_\epsilon, b \Omega^\epsilon$, respectively. In addition, $w_\epsilon, w^\epsilon \rightarrow 1$ a.e. and in every $L^q(b\Omega),$ $1 < q < \infty.$ 
\vskip0.1in
\item If $T_\epsilon, T^\epsilon$ denote the unit tangent vectors for $b \Omega_\epsilon, b \Omega^\epsilon$, respectively, then $T_\epsilon(\Lambda_\epsilon( \cdot))$ and $T^\epsilon(\Lambda^\epsilon(\cdot))$ converge a.e. and in every $L^q(b\Omega)$, $1 < q <
\infty$, to the unit tangent vector for $b \Omega$, $T(\cdot).$

\end{enumerate}
Moreover, it can be deduced from the proof of the construction in \cite{VerchotaThesis}, using the continuity of $\chi_\epsilon$ in the parameter $\epsilon$, that there exists $\delta>0$ small so that for each $y \in \Omega$ with $\operatorname{dist}(y, b \Omega)<\delta,$ then there exists $\epsilon_y>0$ and $y' \in b \Omega$ so that

\begin{equation} y= \Lambda_{\epsilon_y}(y'). \label{Exhausting Curves}\end{equation}
\color{black}

\subsection{Hardy spaces for $\Om$ and for the interior of its complement; the Cauchy Transform for $\Om$} Let $\Omega$ be a Lipschitz domain bounded by a Jordan curve, and let $u$ be defined on $\Omega$ (or $\mathbb{C} \setminus \overline{\Omega}$). We denote by $u^*_i$ and $u^*_e$ the  interior and exterior non-tangential maximal functions of $u$, namely:

$$u^*_i(p):= \sup_{z \in \Gamma_i(p)} |u(z)|, \quad p \in b \Omega; $$ 

$$u^*_e(p):= \sup_{z \in \Gamma_e(p)} |u(z)|, \quad p \in b \Omega $$

where $\Gamma_i(p)$ and $\Gamma_e(p)$ are as in Section \ref{ssec3.1}.

Given $1<p<\infty$, the Hardy space $H^p( \Omega)$ may be defined in either of the following three equivalent ways. Let $f$ be a holomorphic function on $\Omega.$ Then the following three conditions on $f$ are known to be equivalent: (see \cite{L2000}):
\begin{enumerate}
\item There exists a function $g \in L^p(b \Omega)$ which satisfies $\mathbf{C}g=f$, where $\mathbf{C}$ denotes the Cauchy integral for $b\Om$;
\vskip0.1in
\item  $(f)^*_i \in L^p(bD)$;
\vskip0.1in
\item $\displaystyle{\sup_{\epsilon>0} \int\limits_{b\Omega_\epsilon} |f|^p \, d\sigma_\epsilon<\infty}$,
\vskip0.1in
\noindent where $\{\Om_\epsilon\}_\epsilon$ is (any) exhaustion of $\Om$ by
 (rectifiable) subdomains, so here we take $\{\Om_\epsilon\}_\epsilon$ to be as in Section \ref{ssec3.1}.
\end{enumerate}
In either case, we say that $f \in H^p(\Omega)$. It is well-known that such an $f$ has non-tangential boundary limits almost everywhere, and conversely the values of $f$ inside $\Omega$ can be recovered by the Cauchy integral over the boundary. The associated boundary function $\dot{f}$ belongs to $L^p(b \Omega)$, and in particular we write $\dot{f} \in h^p(b \Omega)$ to denote that the function belongs to the boundary Hardy space.  By an important theorem of Smirnov, see 
\cite[Theorem 10.3]{DU}, there is a bijective correspondence between $H^p(\Omega)$ and $h^p(b\Omega)$: on account of this and of condition (1) above, we define the Cauchy transform associated to
$\Om$ as
$$
\mathcal C f(w) \ :=\ (\dot{\mathbf C f }) (w),\quad \text{a. e.}\ w\in b\Om.
$$
It was proved in \cite{CMM} that $\mathcal C$ is bounded: $L^p(b\Om)\to L^p(b\Om)$ for $1<p<\infty$.

Exterior Hardy spaces, which we denote $H^p_{-}( \Omega)$, can be studied in a similar way on $\mathbb C \setminus \overline{\Omega}$. Define the exterior Cauchy integral:

$$ \mathbf{C}_e f(z):= -\frac{1}{2 \pi i} \int\limits_{b \Omega} \frac{f(w)}{w-z} \, dw, \quad z \in \mathbb{C} \setminus \overline{\Omega},$$ where the boundary integration is understood to be in the counterclockwise direction. 
It is immediate that  if $f \in L^1(bD)$, then $\mathbf{C}_ef$ is holomorphic on $\mathbb{C} \setminus \overline{D}$ and vanishes at $\infty.$
In fact, we have the following set of equivalences for a holomorphic function $f$ on $\mathbb{C} \setminus \overline{\Omega}$ which vanishes at infinity: 
\begin{enumerate}
\item There exists a function $g \in L^p(b \Omega)$ which satisfies $\mathbf{C}_e g=f$, where $\mathbf{C}_e$ denotes the exterior Cauchy integral;
\vskip0.1in
\item  $(f)^*_e \in L^p(bD)$;
\vskip0.1in
\item $\displaystyle{\sup_{\epsilon>0} \int\limits_{b\Omega^\epsilon} |f|^p \, d\sigma^\epsilon<\infty,}$\quad with $\{\Omega^\epsilon\}_\epsilon$ as in Section \ref{ssec3.1}.
\end{enumerate}

The proof of \cite[Lemma 2.12]{L2000} shows that the condition $f=\mathbf{C}_e g$ is equivalent to the condition $(f)_e^* \in L^p(b \Omega)$, i.e. the exterior non-tangential maximal function belongs to $L^p(\Omega).$ It is obvious that if said maximal function condition holds, then $\sup_{\epsilon>0} \int\limits_{b\Omega^\epsilon} |f|^p \, d\sigma_\epsilon<\infty.$ Finally, the exhausting domain condition $(3)$ implies the Cauchy integral representation $(1)$ by \cite[Theorem 10.4]{DU} along with an inversion that maps $\mathbb{C} \setminus \overline{\Omega}$ to the interior of a simply connected domain bounded by a Jordan curve.  

Finally, we make a small remark that the definitions here given for $h^p(b\Om)$ and  $h^p_{-}(b \Om)$ can be extended to the endpoint $p=1$ without issue. The equivalence of conditions (1) and (2) is now provided by a result of Jerison and Kenig in \cite{JK}, it is clear (2) still implies (3), and (3) implies (1) again by \cite[Theorem 10.4]{DU}. In fact, conditions (1) and (2) remain equivalent for a broader class of locally rectifiable domains known as chord-arc domains.

\subsection{Proof of Theorem \ref{T:3} and of Theorem \ref{T:1}: Part {\tt i}.}

 In what follows below, we focus on the case $p=2$, but the argument works basically the same for $1<p<\infty.$ The only modification is to replace the Cauchy-Schwarz inequality by H\"{o}lder's inequality, and to use the $L^p$ bound for the Cauchy transform. There are some small modifications to obtain the uniqueness result in the case $p=1$, which we later point out. 

Take $u \in L^2(b\Omega)$. First, to show that such a decomposition exists, we can simply write $u= \mathcal{C}u +(I-\mathcal{C})u$, where $\mathcal{C}$ denotes the Cauchy transform. It is straightforward to check that $(I-\mathcal{C})u$ belongs to $h^2_-(b \Omega)$, while $\mathcal{C} u \in h^2(b \Omega)$ by definition.

We turn now to the uniqueness of this decomposition. Fix $h \in h^2(b \Omega)$ and a Lipschitz boundary parametrization $z(t), a \leq t \leq b$ with counterclockwise orientation, and let $z_0=z(a) \in b\Omega$ denote a point of the boundary where $h$ has a non-tangential boundary limit as before. We will give a slightly different argument here which avoids the construction of curvilinear rectangles (which does not naturally generalize), but instead uses the approximating domains $\Omega_\epsilon$ and continuous family of (smooth) boundary parametrizations for $\Omega_\epsilon$, which we define  $$z_\epsilon(t)=\Lambda_\epsilon(z(t))), \quad a \leq t \leq b. $$

We now turn to the construction of the continuous antiderivative $H(z).$ On the boundary, $H$ is still given by 

$$ H(z(\tau))= \int\limits_{a}^{\tau} h(z(t)) z'(t) \, dt.$$ The Cauchy-Schwarz inequality again yields $H$ is continuous on $b \Omega$. In fact, something stronger is true; since $H$ can be written as an integral of an $L^1$ function with respect to arc length $ds,$ we get that $H$  is absolutely continuous with respect to arc length. To extend the definition of $H$ to $\Omega$, fix $z \in \Omega$ and choose $\varepsilon>0$ sufficiently small so $z \in \Omega_\varepsilon$. Let $\gamma$ be a simple curve contained in $\Omega$ connecting $\Lambda_\varepsilon(z_0)$ and $z$. Set 
$$H(z)= \int\limits_{\eta(z,\epsilon,\gamma)} h(w) \, dw,$$
where $\eta(z, \epsilon, \gamma)$ is the path that starts at $z_0$, travels in a straight line to $\Lambda_\epsilon(z_0)$, and then follows the curve $\gamma$ from the terminus of the short line segment to the point $z$. By convexity, $\eta(z, \epsilon, \gamma)$ always remains within $\Gamma_i(z_0)$ in the first line segment. Note that $$\sup_{w \in \Gamma_i(z_0)}|h(w)|<\infty$$ by virtue of the fact that $h$ has a non-tangential limit at $z_0$, so the integral defining $H(z)$ converges absolutely. It is straightforward to check, using the convexity of $\Gamma_i(z_0)$ and the non-tangential convergence of $h$ at $z_0$, that this definition of $H(z)$ is independent of the choice of $\epsilon$ and of the interior curve $\gamma$. Indeed, given small $\epsilon_1, \epsilon_2$ and arbitrary paths $\gamma_j$ that connect $\Lambda_{\epsilon_j}(z_0)$ and $z_0$, respectively, we have

$$ \int\limits_{\eta(z, \epsilon_1, \gamma_1)} h(z) \, dz-  \int\limits_{\eta(z, \epsilon_2, \gamma_2)} h(z) \, dz = \int\limits_{\mathcal{B}(z,\epsilon_1, \epsilon_2, \gamma_1, \gamma_2)} h(z) \, dz + \int\limits_{\mathcal{T}(z,\epsilon_1,\epsilon_2)} h(z) \, dz  ,$$ 
where $\mathcal{B}(z,\epsilon_1, \epsilon_2, \gamma_1, \gamma_2)$ is the closed curve (not necessarily simple) which consists of $\gamma_1$, followed by $\gamma_2$ traversed backwards, followed by the straight line segment that connects $\Lambda_{\epsilon_2}(z_0)$ to $\Lambda_{\epsilon_1}(z_0)$, and $T(z,\epsilon_1, \epsilon_2)$ is the triangle with ordered vertices $z_0, \Lambda_{\epsilon_1}(z_0)$, and $\Lambda_{\epsilon_2}(z_0)$. 

The first integral on the closed curve $\mathcal{B}(z,\epsilon_1, \epsilon_2, \gamma_1, \gamma_2)$ is clearly $0$ by convexity and Cauchy's theorem. The fact that $\int\limits_{\mathcal{T}(z,\epsilon_1,\epsilon_2)} h(z) \, dz=0$ can be deduced using a limiting argument, the convexity of $\Gamma_i(p)$, and the local boundedness of $h$ within $\Gamma_i(p).$ It is also straightforward to check that $H(z)$ is complex differentiable on $\Omega$ and $H'(z)=h(z)$. Finally, it is clear by construction that $\lim_{\epsilon \rightarrow 0^{+}} H(\Lambda_\epsilon(z_0))=H(z_0)=0.$

It remains to show that $H$ is uniformly continuous up to the boundary. First, we claim that is sufficient to prove that
\begin{equation} H(\Lambda_\epsilon(z(t))\rightarrow H(z(t))
\quad  \text{uniformly in}\ t \in [a,b]\quad \text{as}\  \epsilon \rightarrow 0^{+}.
\label{FlowConvergence} \end{equation} 
 
 A straightforward geometric argument shows that there exists $\delta>0$ so that each cone $\Gamma_i(p)$ is contained inside the non-tangential approach region:
$$\{w\in \Omega: |w-p| \leq (1+\delta)\operatorname{dist}(w, b\Omega)\}. $$ Now consider a general sequence $z_n \rightarrow z \in \partial \Omega$ and let $\varepsilon>0.$ Assuming $n$ sufficiently large, by \eqref{Exhausting Curves}, for each $z_n$, there exists a corresponding point $z_n' \in b \Omega$ and small number $\epsilon_n$ such that $z_n = \Lambda_{\epsilon_n}(z_n')$.  Moreover, by definition and the triangle inequality, $|z_n-z_n'| \leq (1+\delta) |z-z_n|$ and $|z-z_n'| \leq  (2+\delta) |z-z_n|,$ and $|z-z_n|$ can be made arbitrarily small for large $n.$ Using the (uniform) boundary continuity of $H(z)$, choose $N_1$ large enough so that for $n>N_1,$ we have $|H(z)-H(z_n')|< \frac{\varepsilon}{2}$.  Note that $N_1$ only depends on the size of $|z_{N_1}-z|$, not the location of $z$. Then, using the claim \eqref{FlowConvergence}, we may choose $N_2 \geq N_1$ so that for $n>N_2$, $|H(z_n)-H(z_n')|<\frac{\varepsilon}{2}$. Similarly $N_2$ only depends on the size of $|z_{N_2}-z|$, not the particular value of $z$. Then for $n>N_2$, by the triangle inequality there holds $|H(z_n)-H(z)|<\varepsilon,$ completing the proof of the claim.
\vskip0.1in

We turn to the proof of \eqref{FlowConvergence}. To see this, write
\begin{align*}
H(\Lambda_\epsilon(z(t)) & = \int\limits_{L_\epsilon} h(z) \, dz+ \int\limits_{\gamma_\epsilon} h(z)\, d(z),
\end{align*}
where $L_\epsilon$ is given by the straight line segment beginning at $z_0$ and terminating at $\Lambda_\epsilon(z_0)$, and $\gamma_\epsilon$ is the portion of $b\Omega_\epsilon$ beginning at $\Lambda_\epsilon(z_0)$ and terminating at $\Lambda_\epsilon(z(t))$.
Notice that $$ \left| \int\limits_{L_\epsilon} h(z) \, dz \right| \leq |\Lambda_\epsilon(z_0)-z_0| \sup_{w \in \Gamma_i(z_0)} |h(w)| \leq C|\Lambda_\epsilon(z_0)-z_0|,$$
so the first integral converges to $0$ as $\epsilon \rightarrow 0^{+}$ by the convergence of $\Lambda_\epsilon(z_0)$ (independently of $t$).

On the other hand, it can be shown that

$$\int\limits_{\gamma_\epsilon} h(z) \, d(z) \rightarrow \int\limits_{a}^{t} h(z(s)) z'(s)\, ds=H(z(t)) $$
uniformly in $t$ as $\epsilon \rightarrow 0^{+}$ using similar arguments to \eqref{COV}. Indeed, we can estimate using the triangle inequality,

\begin{align*}
\left| H(z(t))-\int\limits_{\gamma_\epsilon} h(z)\, dz  \right| & = \left| \int\limits_{a}^{t} h(z(s)) z'(s) \, ds- \int\limits_{a}^{t} h(z_\epsilon(s)) z_\epsilon'(s) \, ds         \right|\\
& \leq \left| \int\limits_{a}^{t} h(z_\epsilon(s))(z_\epsilon'(s)-z'(s)) \,ds \right| + \left| \int\limits_{a}^{t} (h(z_\epsilon(s))-h(z(s)) z'(s)\,ds \right|\\
& \leq \int\limits_a^b |h(z_\epsilon(s))| |z_\epsilon'(s)-z'(s)| \,ds + \int\limits_a^b |h(z_\epsilon(s))-h(z(s))| |z'(s)| \,ds.
\end{align*}

Note that the right hand side is independent of $t$, so if we prove it tends to $0$ as $\epsilon \rightarrow 0$, we are done. The first integral can be bounded by 

\begin{align*}  \int\limits_{a}^b |h_i^*(z(s))| |z_\epsilon'(s)-z'(s)| \, ds & \leq \|h_i^*\|_{L^2(b \Omega)} \left(\int\limits_{a}^b \frac{| z_\epsilon'(s)-z'(s)|^2}{|z'(s)|^2}|z'(s)| \,ds\right)^{1/2}\\
& \leq   \|h_i^*\|_{L^2(b \Omega)} \bigg( \left(\int\limits_{a}^b |T_\epsilon(z_\epsilon(s))-T(z(s))|^2|z'(s)| \,ds\right)^{1/2}\\
&  + \left(\int\limits_{a}^b \left|\frac{z_\epsilon'(s)}{|z_\epsilon'(s)|}-\frac{z_\epsilon'(s)}{|z'(s)|} \right|^2 |z'(s)|  \,ds\right)^{1/2} \bigg), 
\end{align*}
which tends to zero by the $L^2$ convergence of the unit tangent vectors and the change of variables formula \eqref{COV} together with the convergence of the functions $w_\epsilon$ in $L^2$. 

On the other hand, the second integrand converges pointwise almost everywhere to $0$ by the non-tangential convergence of $h$ at the boundary, and the inequality $$ |h(z_\epsilon(s))-h(z(s))| \leq 2 |h_i^*(z(s))|,$$ which holds for almost every $s$. This observation, together with the Dominated Convergence Theorem, leads to the desired conclusion. 

We make a short remark here to point out that in the case $p=1$, the first integral which is controlled using Cauchy-Schwarz above can instead be bounded directly using only the triangle inequality and the Dominated Convergence Theorem. We leave the details to the interested reader. The second integral may be bounded identically as above. 

Finally, the final conclusion about the uniqueness of the decomposition can be proved similarly to the case when $\Om$ is $C^\infty$ smooth. In particular, we use a theorem, due to Besicovitch,  that states if a function $h$ is holomorphic on $\mathbb{C} \setminus b \Omega$ for $\Omega$ with rectifiable boundary and continuous on $\mathbb{C}$, then $h$ is entire (see \cite{Younsi}). This observation enables us to glue together the interior and exterior holomorphic functions to obtain an entire function that vanishes at infinity, and hence must be identically $0$ by Louiville's Theorem.

\subsection{Proof of Theorem \ref{T:1}: Part {\tt ii.}}

In what follows, we say that a function $u$ defined on $b \Omega$ belongs to the class $C^{0,\alpha}(b\Omega)$ for $0<\alpha < 1$ if $u \circ \gamma \in C^{0,\alpha}[0,1]$, where $\gamma$ is any Lipschitz (that is 
$C^{0, 1}$) parametrization for the boundary of the Lipschitz domain $\Om$.

To prove Theorem {\tt ii.}, it is enough to show that $\mathbf{C}u$ and $\mathbf{C}_e u$ both belong to the class $C^{0,\alpha}(\overline{\Om})$, since we have already proven the uniqueness of the decomposition in $L^2.$ This fact about the regularity of the Cauchy transform is proven in \cite[Theorem 3.3]{Mc}.

\subsection{Proof of Theorem \ref{T:1}: Part {\tt iii.}}

Recall that $\Om$ is said to be of class $C^{k,\alpha}$ for $k\geq 1$ and $0<\alpha < 1$ if there exists a $C^{k, \alpha}$ parametrization of the boundary Jordan curve  $\gamma:[0,1] \rightarrow b \Omega$ such that $\gamma'$ is non-vanishing (and $\gamma$ is of class $C^k$ and moreover 
$\gamma^{(k)}$ is of class $C^\alpha$). 

As before, to prove {\tt iii.}, it suffices to show $\mathbf{C}u \in C^{k,\alpha}(\overline{\Om})$ and $\mathbf{C}_e u \in C^{k,\alpha} (\mathbb{C} \setminus \Om).$ We only prove the first statement, as the argument for $\mathbf{C}_e u$ is the same. Fix $0<\alpha < 1.$  We use will induction on the smoothness index $k$. To establish the base case when $k=1$, we argue as follows. Let $\Om$ be of class $C^{1,\alpha}$ and $u \in C^{1,\alpha}(b \Omega).$ Since $\mathbf{C}u$ is holomorphic on $\Omega$, it is enough to prove $\frac{\partial \mathbf{C}u}{\partial z} \in C^{0,\alpha}(\overline{\Omega})$ (the smoothness up to the boundary is the key part to verify). For fixed $z \in \Omega$, we readily verify by differentiating under the integral sign, parametrizing, and then integrating by parts (the boundary term vanishes since $\gamma$ is closed):

\begin{align*}
\frac{\partial \mathbf{C} u}{\partial z}(z) & = \frac{1}{2 \pi i} \oint_{b \Omega} \frac{u(w)}{(w-z)^2} \,dw   \\
& = \frac{1}{2 \pi i} \int\limits_0^{1} \frac{u(\gamma(t))}{(\gamma(t)-z)^2} \gamma'(t) \, dt\\
&= - \frac{1}{2\pi i} \int\limits_{0}^{1} u(\gamma(t)) \frac{d}{dt}\left(\frac{1}{\gamma(t)-z}\right) \, dt\\
& =\frac{1}{2 \pi i} \int\limits_{0}^1 \frac{g(t)}{(\gamma(t)-z)} \gamma'(t)\, dt\\
& = \mathbf{C}(\widetilde{g})(z)
\end{align*}
where $g(t):=\frac{\frac{d}{dt}\left(u(\gamma(t))\right)}{\gamma'(t)}$ and $\widetilde{u}(\zeta):= g(\gamma^{-1}(\zeta)).$
By definition, it is easy to check $\widetilde{g} \in C^{0, \alpha}(b \Om).$ Since  $\Om$ is Lipschitz by virtue of being class $C^{1,\alpha}$, Theorem \ref{T:1}, part {\tt ii.} then implies $\mathbf{C}\widetilde{g} \in C^{0,\alpha}(\overline{\Omega})$. This observation establishes $\mathbf{C}u \in C^{1,\alpha}(\overline{\Omega}),$ as desired.

Now suppose that the conclusion holds for all $C^{k,\alpha}$ domains and functions for some $k \in \mathbb{N}$, and consider $\Om$ of class $C^{k+1,\alpha}$ and $u \in C^{k+1,\alpha}(b \Omega).$ As in the base case, it suffices to check $\frac{\partial \mathbf{C}u}{\partial z} \in C^{k,\alpha}(\overline{\Omega})$. But we can still obtain the same formula $\frac{\partial \mathbf{C} u}{\partial z}(z)= \mathbf{C}(\widetilde{g})(z),$ where $\widetilde{g}$ has the same definition as above. But this time, we observe $\widetilde{g} \in C^{k, \alpha}(b \Om).$ The induction hypothesis then implies $\frac{\partial \mathbf{C} u}{\partial z} \in C^{k,\alpha}(\overline{\Om})$, completing the proof. 
\section{Further Results}
\label{sec4}
\subsection{Something about Dirichlet and Poisson.} 
As mentioned in the paper \cite{BR} of the same title as this section,
solving the Dirichlet problem and computing the Poisson kernel are
part of the family business, and so it is a pleasant surprise to see
that the $u=h+H$ decomposition provides insight into that old subject,
even in the case of the unit disc. The reader may view this section as an
addendum to \cite{BR}.
\begin{thm}\label{T:4}
Let $u$ be in $C^\infty$ of the unit circle, and let $h$ and $H$ be the terms in the $h+H$ decomposition of $u$. Then,
the solution $U$ of the classical Dirichlet problem for harmonic functions on $\mathbb D$ with datum $u$ admits the following representation:
$$
U(z)=h(z)+H(1/\bar z),\quad z\in \mathbb D.
$$
\end{thm}
Indeed, if $u$ is a $C^\infty$ function on the unit circle, then we have
seen that $u=h+H$ where $h=\mathcal Cu$ is holomorphic inside the
unit disc and $H$ is holomorphic outside the unit disc and tends
to zero at infinity, and both $h$ and $H$ extend smoothly to the
unit circle from their respective sides. The Schwarz function for
the unit disc is $S(z)=1/z$.
Noting that $z=1/\bar z$ on the boundary yields that
$$u(z)=h(z)+H(1/\bar z)$$
when $z$ is on the circle, and this gives us
the harmonic extension $\mathcal E u$ of $u$ to the unit disc
in the form
$$\mathcal E u(z)=h(z)+H(1/\bar z),$$
since $H(1/\bar z)$ is an antiholomorphic function on
the unit disc that goes to zero at the origin, and so has a
removable singularity there.
Note that the value of $\mathcal E u$ at the origin is
$(\mathcal C u)(0)=u_0$, the average of $u$ on the unit circle.
If we now insist that $u$ be real valued, we see that $h(z)$ and
$\overline{H(1/\bar z)}$ are holomorphic functions with the
same imaginary part, and so they differ by a constant, which is
easily seen to be equal to $u_0$. We have proved that
$$\mathcal E u = \left(2\text{Re }\mathcal C u\right)-u_0,$$
and so it is possible to solve the Dirichlet problem in a
rather direct way using
the Cauchy transform on the unit disc! We can also deduce
that, if $u$ is merely continuous and real valued, then the
real part of the Cauchy transform of $u$ on the unit disc is
continuous up to the boundary. Furthermore, the Poisson kernel
is given by
$$P(z,e^{it})=\frac{1}{\pi}
\text{Re }\left[\frac{e^{it}}{e^{it}-z} - \frac{1}{2}\right],$$
which is easily seen to agree with the famous formula in
every complex analysis textbook.
\bigskip

{\bf Remark 2}. We take a moment here to consider the problem
of decomposing a function on the boundary of a simply
connected domain $\Om$ with $C^\infty$ smooth boundary into
the sum of two functions, one given by the boundary values of
a function that is harmonic inside $\Om$ and one given by the
boundary values of a function that is harmonic outside. Given
a $C^\infty$ smooth real valued function $u$ on $b\Om$, we
can take the first function to be the Poisson extension of $u$
to $\Om$ and the second function to be zero outside $\Om$.
This makes it clear that our problem is only interesting
from the point of view of finding a {\it natural\/} or
canonical decomposition. Indeed, note that we can solve the
Dirichlet problem on the outside of $\Om$ by using an inversion
$1/(z-a)$ where $a$ is any point in $\Om$, solving the
problem on the resulting bounded domain, then composing with
the inverse of the inversion. Hence, we can always write $u$ as a
random smooth sum $u_1+u_2$ on the boundary, take the Poisson
extension of $u_1$ inside $\Om$ and the Poisson extension of
$u_2$ outside $\Om$ to make our problem seem silly. We have
shown that a not so silly decomposition is obtained by taking
the first function to be the real part of the Cauchy transform
of $u$ on the inside and the real part of $H$ on the outside,
where $u=h+H$ is the decomposition we have been studying above.
We now consider some other interesting possibilities,
but first we describe a $C^\infty$ version of the
Cauchy-Kovalevskaya theorem that follows from our previous work.

Note that by Theorem~\ref{T:2}, we can find functions $\phi$
in $C^\infty(\C)$ that vanish on $b\Om$ such that
$\dee\phi/\dee\bar z$ agrees with a given $C^\infty$ function
to infinite order along $b\Om$. Hence, by taking conjugates,
we can do the same thing with $\dee/\dee z$ in place of
$\dee\phi/\dee\bar z$. Since the Laplacian is equal to $4$ times
$\dee^2/\dee z\dee\bar z$, we can find real valued
functions $\phi$ such that both $\phi$ and its normal derivative
$\dee\phi/\dee n$ vanish on $b\Om$ and $\Delta\phi$ agrees with
a given real valued $C^\infty$
function to infinite order along $b\Om$. Indeed, given
a smooth real valued function $U$, first get a complex
valued function $v$ that vanishes on $b\Om$
such that $\dee v/\dee z$ agrees with $\frac14 U$ to infinite order
along $b\Om$. Next, find a complex valued $\phi$ that vanishes along
$b\Om$ such that $\dee\phi/\dee\bar z$ agrees with
$v$ to infinite order along $b\Om$. Now, it is a rather
straightforward exercise to show that the real
part of $\phi$ is the solution we seek. Indeed, let $\rho$
be the defining function for $\Om$ that is equal to minus
the distance to the boundary on $\Om$ and plus the distance
outside. Since $\phi$ is zero on the boundary, the tangential
derivative $\nabla\phi\cdot(\rho_y,-\rho_x)=\frac{2}{i}
(\phi_z\rho_{\bar z}-\phi_{\bar z}\rho_z)$ is zero along
$b\Om$. Consequently,
$\phi_z\rho_{\bar z}=\phi_{\bar z}\rho_z=v\rho_z=0$
along $b\Om$. The
normal derivative of $\phi$ is $\nabla\phi\cdot\nabla\rho=2
(\phi_z\rho_{\bar z}+\phi_{\bar z}\rho_z)$, which is now seen
to also be zero along $b\Om$. Finally, since $\Delta \phi-U$
vanishes to infinite order along $b\Om$, we see that so
does the real part. Furthermore, the real part of $\phi$ and
its normal derivative are also zero along $b\Om$ and our
claim is proved. We will now apply this result to Green's
identity in a manner similar to the way we used
Theorem~\ref{T:2} in the Cauchy-Pompieu formula.

Green's formula is
$$U(z)=\int_{b\Om}u(w)\frac{\dee}{\dee n_w}N(z,w)-
N(z,w)
\frac{\dee}{\dee n_w}u(w)\ ds_w +
\iint_{\Om}N(z,w)\Delta U(w)\ dA_w,$$
where $\dee/\dee n_w$ represents the outward normal derivative
in the $w$ variable, $dA_w$ is Lebesgue area measure in the
$w$ variable, $ds_w$ is the element of arc length in the $w$
variable, and
$N(z,w)=\frac{1}{2\pi}\text{Ln\,}|z-w|$ is the fundamental
solution for the Laplacian. Given a function $u$ that is
$C^\infty$ smooth on $b\Om$, we can extend $u$ to a function
$U$ in $C^\infty(\Obar)$ in many ways. One way that leads to
interesting results is to extend so that $U=u$ on $b\Om$ and
$\dee U/\dee n$ is zero. (To do this, note that functions of
the form $v(z)\rho(z)$, where $\rho$ is the defining function
mentioned above, are zero on the boundary and have normal
derivative equal to $v$ on the boundary, and so can be used
to construct functions locally with prescribed values and
normal derivatives on the boundary.)
Assuming we have such an extension $U$, we
now choose a function $\phi$ in $C^\infty(\C)$ such that
$\phi$ and its normal derivative vanish on the boundary and
$\Delta\phi$ agrees with $\Delta U$ to infinite order along
$b\Om$. If we apply Green's identity to $U-\phi$, we can deduce
that the double integral extends $C^\infty$ smoothly to $\C$ and
has Laplacian equal to zero outside $\Om$ and equal to
$\Delta(U-\phi)$ in $\Om$.  Hence, we deduce that $u= u_1+u_2$
on the boundary where
$$u_1(z)=
\int_{b\Om}u(w)\frac{\dee}{\dee n_w}N(z,w)\ ds_w$$
is harmonic on $\Om$ and has $C^\infty$ smooth boundary
values and
$$u_2(z)=
\iint_{\Om}N(z,w)\Delta [U(w)-\phi(w)]\ dA_w$$
is harmonic on the outside of $\Om$ and has $C^\infty$ smooth
boundary values along $b\Om$. Note that another of Green's
identities yields that
$$u_2(z)=
\iint_{\Om}N(z,w)\Delta U(w)\ dA_w=
-\int_{b\Om}u(w)\frac{\dee}{\dee n_w}N(z,w)\ ds_w.$$
outside $\Om$. Before anyone gets too excited, we must observe
that the normal derivative of $N(z,w)$ is (see \cite[p.~94]{Be}
for details about this type of calculation)
$$\frac{\dee N}{\dee n}=-i \frac{\dee N}{\dee z} T+
i\frac{\dee N}{\dee\bar z}\,\overline{T}
=2\text{Re}\left(-i \frac{\dee N}{\dee z} T\right),$$
and so
$$\frac{\dee N}{\dee n_w}(z,w)\ ds_w=
\text{Re}\left(\frac{1}{2\pi i} \frac{T(w)}{w- z}\ dw\right),$$
and we see that $u_1$ is just the real part of the Cauchy
transform of $u$ inside $\Om$, and $u_2$ is minus the real
part of the Cauchy transform when $z$ is outside $\Om$.
Thus, this game just yields another way to look at the higher
order version of the Plemelj formula we have mentioned and
alludes to the possibility that the real part of the Cauchy
transform applied to real valued functions might have
special features, as it does in the unit disc.

Another interesting consequence of these ideas is obtained if
we construct the function
$U$ so that $U$ is equal to zero on $b\Om$ and $\dee U/\dee n$
is equal to $u$ on the boundary, then choose $\phi$ in the same
way and restrict Green's identity to the boundary to deduce that
$$\int_{b\Om} N(z,w) u(w)\ ds_w$$
is harmonic on $\Om$ and has $C^\infty$ boundary values that
agree with the $C^\infty$ boundary values of
$$\iint_{\Om}N(z,w)\Delta [U(w)-\phi(w)]\ dA_w=
\iint_{\Om}N(z,w)\Delta U(w)\ dA_w=
\int_{b\Om} N(z,w) u(w)\ ds_w$$
from the outside of $\Om$. (That these functions agree on the
boundary is no surprise since the kernel has an $L^1$ singularity
along the boundary.)

These ideas beg for further study, but we leave this to the
future.
\bigskip

{\bf Remark 3}. If $\Om$ is a simply connected domain with real analytic
boundary, then there is a Schwarz function $S(z)$ that is
analytic in a neighborhood of the boundary such that
$S(z)=\bar z$ on the boundary, and $\overline{S(z)}$ is
an anitholomorphic reflection function of the boundary.
It maps an open subset on the inside of the boundary
curve to an open subset on the outside of the boundary and
it fixes the boundary. When we write
$$u=\mathcal C u + H(\,\overline{S(z)}\,)$$
for $z$ in the boundary in this case, we see that we obtain
a function that has boundary values given by $u$ and that
extends inside $\Om$ as a function that
is harmonic on the open subset of the inside of
the boundary curve where $\overline{S(z)}$ is
antiholomorphic and maps to the outside of the boundary.
Only in the case of a disc can one fill in the whole inside
this way to obtain the harmonic extension to $\Om$.
\bigskip
\subsection{A pseudo-local property of the Szeg\H o projection.}
Because the tools and notation set up in \S\ref{sec2} can
be used to prove a pseudo-local property of the Szeg\H o
projection that has eluded us in the past, we take this
opportunity to describe and prove it. 

\begin{thm}\label{T:5}
Suppose $\Om$ is a domain bounded by a $C^\infty$ smooth Jordan curve.

Let $z_0\in b\Om$, and let $\mathbb D_r(z_0)$ be 
a (small) disc of radius $r$ about $z_0$ such that $\mathbb D_r(z_0)\cap\Om$ is simply connected. 
Let $h$ be holomorphic on $\mathbb D_r(z_0)\cap\Om$
and with $L^2$ boundary limits on $\mathbb D_r(z_0)\cap b\Om$.

Then, there is $f$ in $h^2(b\Om)$ such that the difference $f-h$ extends $C^\infty$ smoothly
to the boundary of $\Om$ near $z_0$.
\end{thm}

If we were thinking about this problem in the Bergman space,
we would use Spencer's formula
$B=I-4\frac{\dee}{\dee z}G\frac{\dee }{\dee\bar z}$
that relates the Bergman projection $B$ to the classical
Green's operator $G$ for the Laplacian with homogeneous
boundary conditions. Pseudo-local properties for the Green's
operator would allow us to use the Bergman projection to
get a function $f$ in the Bergman space that has the same
behavior near $z_0$ modulo holomorphic functions that
extend smoothly to the boundary near $z_0$. These are
standard techniques of PDE. In the Hardy space setting, we
resort to the following methods.

Let $P$ denote the Szeg\H o projection, $S(z,w)$ the
Szeg\H o kernel, and $L(z,w)$ the Garabedian kernel.
Basic properties of these objects are described and
proved in \cite[pp.~13-32]{Be}.
Let $u$ be the function on the boundary of $\Om$
gotten by extending the boundary values of $h$ in
$\mathbb D_r(z_0)\cap b\Om$ to the whole boundary via extension
by zero and let $f=Pu$. Then
$$(Pu)(z)=\int_{w\in b\Om}S(z,w)u(w)\ ds$$
for $z$ in $\Om$.
Pick a point $w_0$ in the clockwise direction from $z_0$ in
$\mathbb D_r(z_0)\cap b\Om$ where $h$
has a non-tangential limit and set up the curves
$\mathcal R(0,\epsilon,\tau)$ as in \S\ref{sec2} using
$w_0$ in place of $z_0$. Pick a $\tau$ so that
$\mathcal R(0,\epsilon,\tau)$ is contained in
$\mathbb D_r(z_0)\cap \Om$ and
the part $\sigma$ of
$\mathcal R(0,\epsilon,\tau)$ in the boundary is
compactly contained in 
$\mathbb D_r(z_0)\cap b\Om$. Let $-\lambda$ denote the
counterclockwise part of
$\mathcal R(0,\epsilon,\tau)$ inside $\Om$ (i.e. the
part with $\sigma$ removed).
We now split the Szeg\H o projection
integral into two pieces and note that
$$S(z,w)\,ds=\frac{1}{i}L(w,z)\,dw$$
to obtain
$$(Pu)(z)=\int_{b\Om-\sigma}S(z,w)u(w)\ ds +
\frac{1}{i}
\int_{\sigma}L(w,z)h(w)\ dw.$$
Assume $z$ is inside $\mathcal R(0,\epsilon,\tau)$.
We know a Cauchy theorem for integrals of holomorphic
functions with $L^2$ boundary values along $b\Om$ for
$\mathcal R(0,\epsilon,\tau)$, and so a residue theorem
for such functions with finitely many poles inside also
holds. The residue of $L(w,z)$ at $w=z$ is $1/(2\pi)$, and so the second
integral is equal to
$$h(z)+
\frac{1}{i}
\int_{\lambda}L(w,z)h(w)\ dw.$$
The formula for $Pu$ now yields that
$$f(z)-h(z)=
\int_{b\Om-\sigma}S(z,w)u(w)\ ds +\frac{1}{i}
\int_{\lambda}L(w,z)h(w)\ dw,$$
and we can read off the fact that this function
extends $C^\infty$ smoothly to the boundary of $\Om$ near
$z_0$ because $S(z,w)$ extends $C^\infty$ smoothly to
$\Obar\times\Obar$ minus the boundary diagonal, and
$L(z,w)$ extends $C^\infty$ smoothly to
$\Obar\times\Obar$ minus the diagonal of $\Obar$ (see
\cite[pp.~141-145]{Be} for proofs of these facts).
The proof is complete.


\begin{thebibliography}{000}

\bibitem{Be}
Bell, S.,
{\em The Cauchy transform, potential theory, and
conformal mapping}, 2nd Edition, CRC Press, Boca Raton, 2016.

\bibitem{Be1}
Bell, S.,
{\em Unique continuation theorems for the $\bar\partial$-operator
and applications}, J. of  Geometric Analysis
{\bf 3} (1993), 195--224.

\bibitem{BR}
Bell, S. and Reyna de la Torre, L., {\em Something about Poisson
and Dirichlet}, Chapter 1 in Handbook of Complex Analysis,
Steven G. Krantz, editor, pp.~1-17, Taylor and Francis, 2022.

\bibitem{Ca} Calder\'{o}n A.,
{\em Cauchy integrals on Lipschitz curves and related operators},
Proc. Nat. Acad. Sci. U.S.A. {\bf 74} no.4 (1977), 1324--1327.

\bibitem{CMM} Coifman, R.R., McIntosh, A., and Meyer, Y. 
{\em L'Int\`egrale de Cauchy d\`efinit un op\`erateur born\`e sur $L^2$  pur les courbes Lipschitziennes}, Ann. Math. \textbf{116} (1982), 361--387.

\bibitem{D} David G.,
{\em Op\`erateurs int\`egraux singuliers sur certain courbes due plan complexe}, Ann. Sci. \'{E}cole Norm. Sup. (4) \textbf{17} (1984), 157--189.

\bibitem{DU} Duren, Peter L. {\em Theory of $H^p$ spaces}.
Pure Appl. Math., Vol. 38 Academic Press, New York-London, 1970.


\bibitem{GG}
Guillemin, V. and Golubitsky, M., {\em Stable mappings and their
singularities}, Springer Graduate Texts in Math., New York, 1973.

\bibitem{HLLS1} Hulse J, Lanzani L., Luca E. and Llewellyn Smith S.,
{\em A Transform-based technique for solving boundary value problems on convex planar domains}, 
IMA J. of Applied mathematics {\bf 89} no. 3 (2024) 574 - 597



\bibitem{HLLS2} Hulse J, Lanzani L., Luca E. and Llewellyn Smith S.,
{\em The unified transform method: beyond circular or convex}, submitted (ArXiv: 2410.09744).

\bibitem{JK} 
Jerison, D.S. and Kenig, C.E. {\em Hardy spaces, $A_\infty$, and singular integrals on chord-arc domains}, Math. Scand. \textbf{50} (1982), no. 2, 221–247.

\bibitem{K}
Koosis, P., {\em Introduction to $H_p$ spaces}, Cambridge University
Press, London Math. Soc. Lecture Note Series 40, Cambridge, 1980.

\bibitem{L2000} Lanzani, L. {\em Cauchy transform and Hardy spaces for rough planar domains}, Contemp. Math. 251, pp. 409-428, American Mathematical Society, Providence, RI, 2000.

\bibitem{LS}
Lanzani, L., and Stein, E. M.,
{\em Szeg\"{o} and Bergman projections on non-smooth planar domains},
J. Geom. Anal. {\bf 14} (2004), 63--86.

\bibitem{LS2} Lanzani,L., and Stein, E.M.,
{\em The Cauchy integral in $\mathbb{C}^n$ for domains with minimal smoothness}, Adv. Math. {\bf 264} (2014), 776–-830.

\bibitem{Mc}
McLean, W. {\em H\"older estimates for the Cauchy Integral on a Lipschitz contour}, J. of Integral Equations and Applications {\bf 1} no. 3  (1988), pp.435 - 451.

\bibitem{P} I. I. Privalov, {\em Boundary Properties of Analytic Functions}, GITTL, Moscow-Leningrad, 1950.

\bibitem{Verchota} Verchota, G., {\em Layer potentials and regularity for the Dirichlet problem for Laplace's equation in Lipschitz domains}, J. Funct. Anal. {\bf 59} (1984), 572--611.

\bibitem{VerchotaThesis} Verchota, G., {\em Layer Potentials and Boundary Value Problems for Laplace's Equation on Lipschitz Domains }, PhD Thesis, pp.~107--120, ProQuest LLC, Ann Arbor, MI, 1982.

\bibitem{Younsi} Younsi, M. \em {On removable sets for holomorphic functions}, EMS Surv. Math. Sci. {\bf 2}(2015), 219--254.

\end{thebibliography}
\end{document}